\documentclass[12pt]{article}
\usepackage{amsmath}
\usepackage{latexsym}
\usepackage{amssymb}
\newcommand{\semid}{\mbox{$\times\!$\rule{.15 mm}{2.07 mm}}}
\newcommand{\dl}{{\displaystyle \lim_{\rightarrow}}}
\newcommand{\pl}{{\displaystyle \lim_{\leftarrow}}}
\newcommand{\sbull}{{\scriptscriptstyle \bullet}}

\begin{document}
\begin{center}
{\bf\Large Fundamental Problems in the Theory of\vspace{1.7mm}
Infinite-Dimensional Lie Groups}\\[4mm]
{\bf Helge Gl\"{o}ckner}
\end{center}
\begin{abstract}
\noindent
In a preprint from 1982,
John Milnor formulated various fundamental
questions concerning infinite-dimensional
Lie groups.
In this note, we describe some of the answers
(and partial answers)
obtained in the
preceding years.\vspace{2mm}
\end{abstract}
\begin{center}
{\bf\Large Introduction}
\end{center}
Specific classes of infinite-dimensional Lie groups
(like groups of operators,
gauge groups, and diffeomorphism groups)
have been studied extensively and
are well understood.
But much less is known about
general infinite-dimensional Lie groups,
and many fundamental problems
are still unsolved.
Typical problems were recorded
in John Milnor's preprint~\cite{PRE},
which preceded his well-known survey article~\cite{Mil}.
In this note, we recall
Milnor's questions
and their background
and describe some of the answers
(or partial answers)
obtained so far.
\section{Basic definitions}
To define infinite-dimensional Lie groups,
John Milnor uses
the following
notion of smooth maps
between locally convex spaces
(which are known as ``Keller $C^\infty_c$-maps'' \cite{Kel}
in the literature):\\[2.5mm]
{\bf Definition.}
Let $E$ and $F$ be real locally convex spaces,
$U\subseteq E$ be open,
and $f\colon U\to F$ be a map.
For $x\in U$ and $y\in E$,
let
$(D_yf)(x):=\frac{d}{dt}\big|_{t=0}\, f(x+ty)$
be the directional derivative
(if it exists).
Given $k\in {\mathbb N} \cup\{\infty\}$,
the map $f$ is called $C^k$
if it is continuous,
the iterated directional derivatives
\[
d^jf(x,y_1,\ldots, y_j):=(D_{y_j}\cdots D_{y_1}f)(x)
\]
exist for all $j\in {\mathbb N}$ such that $j\leq k$,
$x\in U$ and $y_1,\ldots, y_j\in E$,
and all of the maps
$d^jf\colon U\times E^j\to F$
are continuous. As usual, $C^\infty$-maps
are also called \emph{smooth}.\\[2.5mm]
A \emph{smooth manifold} modeled on a locally convex topological
vector space~$E$ is
a Hausdorff topological space~$M$,
together with a set ${\mathcal A}$
of homeomorphisms from
open subsets of~$M$ onto open subsets
of~$E$, such that the domains cover~$M$
and the transition maps are smooth.
Smoothness of maps between manifolds
is defined as in the finite-dimensional
case (it can be tested in local charts).
Also products of manifolds are defined as usual.
A \emph{Lie group} is a group, equipped with
a smooth manifold structure modelled
on a locally convex space~$E$
such that the group
operations are smooth maps.
Lie groups modelled on Banach spaces
are called \emph{Banach-Lie groups}.
As in finite dimensions,
the tangent space
$L(G):=T_1(G)\cong E$ at the identity element
of a Lie group~$G$
can be made a topological Lie algebra
via the identification with the Lie algebra
of left invariant vector fields on~$G$.\\[2.5mm]
Milnor requires in \cite{PRE}
that the modelling space
is complete, and relaxes the condition
to sequential completeness (convergence of
Cauchy sequences) in \cite{Mil}.
We follow his custom here
(unless we explicitly state the contrary).\\[2.5mm]
Occasionally, we shall
also encounter
analytic mappings
and the corresponding Lie groups.
Given complex locally convex spaces
$E$ and~$F$,
a map $f\colon U\to F$
on an open subset $U\subseteq E$ is called
\emph{complex analytic}
if it is continuous
and for each $x\in U$
there exist a $0$-neighbourhood
$Y\subseteq E$ and
continuous homogeneous polynomials
$\beta_n\colon E\to F$ of degree~$n$
such that $x+Y\subseteq U$ and
\[
f(x+y)\;=\; \sum_{n=0}^\infty\, \beta_n(y)\quad
\mbox{for all $y\in Y$}
\]
as a pointwise limit
(see \cite{BaS}).
Given real locally convex spaces,
following \cite{Mil}
a map $f\colon U\to F$
on an open subset $U\subseteq E$
is called \emph{real analytic} if it extends
to a complex analytic map
between open subsets of the complexifications
$E_{\mathbb C}$ and $F_{\mathbb C}$.
We remark that the above definition
of $C^k$-maps also makes sense
over the complex field of scalars;
it is known that mappings to sequentially complete
complex locally convex spaces are complex analytic
if and only if they are $C^1$ in the complex sense.
Further information can be found in~\cite{GN1}.
In particular, complex analytic maps
are real analytic, and
real analytic maps are smooth.
\section{Existence of an exponential map}
Let $G$ be a Lie group.
Given $X\in L(G)$,
there is at most one smooth homomorphism
$\gamma_X\colon {\mathbb R}\to G$
with $\gamma_X'(0)=X$.
If $\gamma_X$ always exists,
$G$ is said to
\emph{have an exponential map}, and we define
it via
$\exp_G\colon L(G)\to G$,
$\exp_G(X):=\gamma_X(1)$.\\[2.5mm]
Milnor asked \cite[p.\,1]{PRE}
whether every Lie group
has a smooth exponential map.
This question is still wide open:
neither is it known whether
an exponential map always exists,
nor whether smoothness is automatic.\\[2.5mm]
In the absence
of completeness properties
of the modelling space of a Lie group
(which Milnor requires),
an exponential mapping need not exist
(see \cite[\S6]{ALG}):\\[2.5mm]
{\bf Example.}
Consider the algebra ${\mathbb R}[X]$ of polynomial
functions $[0,1]\to{\mathbb R}$
and the algebra of fractions
$A:=S^{-1}{\mathbb R}[X]\subseteq  C[0,1]$,
where $S$ is the set of all polynomial
functions without zeros in~$[0,1]$.
Then $A$ is a non-complete
topological algebra in the topology
induced by the Banach algebra $C[0,1]$.
Since $A^\times=A\cap C[0,1]^\times$,
the unit group
$A^\times$ is open in~$A$.
Hence $A^\times$ is a Lie group.
It does not have an exponential
map since $\gamma_f$ only exists for
$f\in {\mathbb R}{\bf 1}\subseteq  A=L(A^\times)$.\\[2.5mm]
Of course, a smooth exponential map
does exist for all typical classes of Lie groups.
But the construction of $\exp_G$
(and its particular properties)
strongly depend on the type of Lie group:\\[2.5mm]
{\bf Banach-Lie groups.}
As a consequence
of the local existence and uniqueness
of solutions to ordinary differential equations
in Banach spaces,
every Banach-Lie group~$G$ has
a smooth exponential map (cf.\ also Section~\ref{secreg}).
Since $T_0(\exp_G)=\text{id}_{L(G)}$,
the inverse function theorem for
smooth maps between Banach spaces
implies that $\exp_G$ is a local diffeomorphism
at~$0$.\\[2.5mm]
{\bf Linear Lie groups.}
Let $A$ be a \emph{continuous
inverse algebra}, viz.\
a locally convex topological
algebra whose unit group $A^\times$ is open and such that
the inversion map
$\iota \colon A^\times \to A$,
$x\mapsto x^{-1}$
is continuous.
Then $\iota$ is analytic and
thus $A^\times$ is an analytic Lie group.
If $A$ is sequentially complete,
then the exponential series converges
and defines an analytic map
$\exp
\colon A\to A^\times$,
$\exp(x):=\sum_{n=0}^\infty
\frac{1}{n!}x^n$
which is the exponential
map of~$A^\times$ (see \cite[Theorem~5.6]{ALG}).
After replacing $A$ with $A_{\mathbb C}$
if necessary,
this follows from the fact that
\[
\exp(x)\,=\, \frac{1}{2\pi i}\int_{|\zeta|=r}
e^{\,\zeta}\cdot (\zeta-x)^{-1}\, d\zeta
\]
for $x\in A$
in terms of holomorphic functional
calculus, where $r$ is chosen so large
that the circle $|\zeta|=r$ surrounds
the spectrum of~$x$.
Here $\exp$ is a local diffeomorphism,
with $\exp^{-1}(x)=
\log(x)=\sum_{n=1}^\infty \frac{(-1)^{n+1}}{n}(x-1)^n$
for $x$ near~$1$.\\[2.5mm]
{\bf Mapping groups.}
Let $M$ be a compact manifold,
$G$ a Lie group
with a smooth exponential map
$\exp_G$ (e.g.,
a finite-dimensional Lie group).
Then
$C^\infty(M,G)$ is a group with pointwise
group operations, and can be made a Lie group
modelled on $C^\infty(M,L(G))$.
The map
$C^\infty(M,L(G)) \to
C^\infty(M,G)$,
$\gamma \mapsto  \exp_G\circ \, \gamma$
is smooth and is easily seen to be the
exponential map of $C^\infty(M,G)$
(cf.\ \cite{Mil}, \cite{GCX}).
If $G$ has a locally diffeomorphic
exponential map, then also $C^\infty(M,G)$.\\[2.5mm]
{\bf Diffeomorphism groups.}
For each compact smooth manifold~$M$,
the group $G:={\rm Diff}(M)$ of $C^\infty$-diffeomorphisms
of~$M$ can be made a Lie group (with composition
as the group multiplication),
modelled on
the space ${\mathcal V}(M)$ of smooth vector fields.
It has a smooth exponential map given by
\[
\exp_G\colon {\mathcal V}(M)\to G\,,\quad
X\mapsto \Phi_X(1,\sbull)\,,
\]
where $\Phi_X\colon {\mathbb R}\times M\to M$
is the flow of vector field~$X$ (see
\cite{Mic} and \cite{Mil}).
Already for $M:={\mathbb S}^1$,
$\exp_G$
is not a local diffeomorphism at~$0$ (see \cite[p.\,1017]{Mil}).\\[2.5mm]
{\bf Direct limit groups.}
Given an ascending sequence
$G_1\subseteq G_2\subseteq\cdots$ of finite-dimensional
Lie groups such that the
inclusion maps are smooth homomorphisms,
we consider $L(G_n)$ as a Lie subalgebra
of $L(G_{n+1})$.
Then
$G:=\bigcup_{n\in {\mathbb N}} G_n=\dl\, G_n$\vspace{-.5mm}
is a group in a natural way, which
can be given a smooth manifold
structure modelled
on the locally convex direct limit
${\mathfrak g}:=\dl\, L(G_n)$\vspace{-.5mm}
making it Lie group
(see \cite{FUN}; cf.\
\cite{NRW} for an earlier, more restricted
method).
The map ${\mathfrak g}\to G$,
$x\mapsto \exp_{G_n}(x)$ if $x\in L(G_n)$
is the exponential map of~$G$;
it is smooth.
\section{Analyticity of multiplication in exponential coordinates}
As illustrated by the examples above,
many (but not all) infinite-dimensional
Lie groups $G$ are \emph{locally exponential}
in the sense that $\exp_G$ exists and is
a local $C^\infty$-diffeomorphism at~$0$.
A locally exponential Lie group~$G$
is called a \emph{BCH-Lie group}
if the group multiplication
is analytic in exponential
coordinates,
i.e., if
$(x,y)\mapsto x*y:=\exp_G^{-1}(\exp_G(x)\exp_G(y))$
is analytic
on some open $0$-neighbourhood in $L(G)\times L(G)$.
Then $x*y$ is given by the Baker-Campbell-Hausdorff
(BCH-) series~\cite{GN2}.
In our terminology,
Milnor asked (cf.\ \cite[p.\,31]{PRE}):\\[2.5mm]
\emph{If {\rm\,(a)\,} $G$ is locally exponential,
or {\rm\,(b)\,} $G$ is real or complex analytic,
does it follow that $G$ is BCH\,}?\\[2.5mm]
The answers to both questions are negative.\\[2.5mm]
(a) A counterexample for (a) is mentioned in \cite[p.\,823]{Rob}.
Slightly simpler is
\[
G := {\mathbb R}^{\mathbb N}\semid_\alpha\,  {\mathbb R}
\quad\mbox{with}\quad
\alpha(t).(x_n)_{n\in {\mathbb N}}:=(e^{nt} x_n)_{n\in {\mathbb N}}.
\]
Using the identity map
onto the Fr\'{e}chet
space ${\mathbb R}^{\mathbb N}\times {\mathbb R}$
as a global chart,
$G$ becomes a real analytic Lie group.
Its exponential map
$\exp_G\colon {\mathbb R}^{\mathbb N}\semid\,{\mathbb R}\to G$,
$\exp_G((x_n)_n,t)
=\big(\big(\frac{e^{nt}-1}{nt} x_n\big)_n, t\big)$
is a $C^\infty$-diffeomorphism
and real analytic
(whence $G$ is locally exponential),
but $\exp_G^{-1}$ is not real analytic.
It can be shown that the map
$(x,y)\mapsto x*y$ is not real analytic
on any $0$-neighbourhood
in $L(G)\times L(G)$, and thus $G$
is not a BCH-Lie group
(cf.\ \cite{GN2} for details).\\[2.5mm]
(b) The group $G:= {\mathbb C}^{({\mathbb N})}\semid_\alpha \, {\mathbb R}$
with $\alpha(t).(x_n)_{n\in {\mathbb N}}
:=(e^{int} x_n)_{n\in {\mathbb N}}$
is a real ana\-lytic Lie group with a global chart,
the identity map onto
${\mathbb C}^{({\mathbb N})}\times {\mathbb R}$
(equipped with the finest locally convex vector
topology).
As shown in \cite[Example~5.5]{DIR},
the exponential map
$\exp_G((x_n)_n,t)
=\big(\big(\frac{e^{int}-1}{int} x_n\big)_n, t\big)$
is not injective on any $0$-neighbourhood,
and the exponential image is not an identity neighbourhood.
Hence $G$ is not~BCH (not even locally exponential).
The corresponding semi\-direct product
${\mathbb C}^{({\mathbb N})}
\semid \,  {\mathbb C}$
provides a complex analytic
counterexample.\\[2.5mm]
The counterexample for (a) was stimulated by
Neeb's discussions of projective limits
of finite-dimensional Lie groups~\cite{GN2}.
Our counterexamples show
that neither projective nor direct
limits of finite-dimensional
(and hence BCH-) Lie groups
need to be locally exponential.
For further information on
BCH-Lie groups and locally exponential
Lie groups, see \cite{GCX}, \cite{ALG}, \cite{GN2},
\cite{PRE}, \cite{Rob} and \cite{Rb2}.
\section{Regularity questions}\label{secreg}
Roughly speaking,
a Lie group $G$ is
called regular
if all ODEs of interest
for Lie theory can be solved in~$G$,
and the solutions depend smoothly on parameters.
Formally,
a Lie group
is \emph{regular}
if the following holds
(see \cite[Definition~7.6]{Mil}):
\begin{itemize}
\item[(a)]
Every smooth curve
$\gamma\colon [0,1]\to L(G)$
arises as the left logarithmic
derivative
of a (necessarily unique)
smooth curve $\eta\colon [0,1]\to G$,
that is, $\gamma(t)=
\eta(t)^{-1}\cdot \eta'(t)$
for all $t\in [0,1]$ (taking the product
in the Lie group~$TG$);
\item[(b)]
The mapping $C^\infty([0,1],L(G))\to G$
taking $\gamma$ to~$\eta(1)$
is smooth, where $C^\infty([0,1],L(G))$
is equipped with its usual
locally convex topology.
\end{itemize}
Regularity is a useful property.
For example, every regular Lie group
has a smooth exponential map.
Also, every continuous
homomorphism
\[
\phi\colon L(G)\to L(H)\,,
\]
where $G$ is a simply connected Lie group
and~$H$ regular, gives rise
to a unique smooth homomorphism $\psi\colon G\to H$
with $T_1\psi=\phi$ (see \cite[Theorem~8.1]{Mil};
cf.\ \cite[Theorem~5.4]{PRE} for a
precursor by Thurston
for so-called ``receptive''
Lie groups, which coincide
with regular Lie groups by \cite[Lemma~8.8]{Mil}).\\[2.5mm]
It is unknown whether every Lie group
is regular, although all typical examples
are regular:
Regularity of Banach-Lie groups
follows from the smooth dependence
of solutions to ODEs in Banach spaces
on parameters;
regularity of ${\rm Diff}(M)$
for compact~$M$
was proved in~\cite{Mil} (see also \cite{OYY},
where an earlier, stronger notion of regularity
was used);
and regularity of $C^\infty(M,G)=\bigcap_{k\in {\mathbb N_0}}C^k(M,G)$
with finite-dimensional~$G$
can be reduced to the Banach case.\\[2.5mm]
Also the group $\text{Diff}_c(M)$
of compactly supported smooth diffeomorphisms
of a $\sigma$-compact finite-dimensional
smooth manifold~$M$ can be made
a Lie group, and in fact in two ways:
It can be modelled either on the LF-space
${\mathcal V}_c(M)=\dl_K {\mathcal V}_K(M)$\vspace{-.4mm}
of compactly supported smooth vector fields,
equipped with the
locally convex direct limit topology;
or on the same vector
space, but equipped with the coarser
topology making it the projective limit
\[
{\mathcal V}_c(M)
\;=\; \bigcap_{k\in {\mathbb N}_0}{\mathcal V}^k_c(M)
\;=\; \pl_{k\in {\mathbb N}_0}{\mathcal V}^k_c(M)
\]
of the LB-spaces of compactly supported
$C^k$-vector fields.
The first discussion of $\text{Diff}_c(M)$
was given in \cite{Mic} (even for paracompact manifolds).
A different, more elementary construction
was described later in~\cite{DIF}.
The regularity of both Lie group structures
on $\text{Diff}_c(M)$ was asserted in~\cite{PRE}
(using other terminology)
and fully proved in~\cite{DIF}.\\[2.5mm]
Also every direct limit group (as described in Section~2)
is regular, by \cite[Theorem~8.1]{FUN}.
The unit groups of
sequentially
complete continuous inverse algebras
are regular
as a consequence of
results by Robart~\cite{Rb2},
who addressed the question whether
every BCH-Lie group
is regular and achieved essential
progress in this direction.\\[2.5mm]
See~\cite{KaM} for a counterpart of
regularity in the convenient setting of analysis.
As in the case of convenient regularity~\cite{MaT},
an abelian
Lie group~$G$ modelled on a Mackey complete
locally convex space~$E$ is regular if and only if
$G\cong E/\Gamma$
for a discrete subgroup $\Gamma\subseteq E$
(see \cite[Proposition~V.1.9]{MON} or \cite{GN2}).
Neeb also showed that every solvable Lie group
with smooth exponential map
is regular~\cite{GN2}.\\[2.5mm]
Criteria for convenient regularity
were given in~\cite{Tei} and
applied to the ``strong ILB-Lie groups''
of Omori and collaborators
(as in \cite{Omo}).\\[2.5mm]
Related to regularity is another question
by Milnor \cite[p.\,1]{PRE}:
{\em If two simply\linebreak
connected Lie groups
$G$ and $H$ have isomorphic Lie algebras,
does it follow that
$G\cong H$}\,?
The theorem by Thurston and Milnor just described
implies that the answer is positive
if both $G$ and $H$ are regular
\cite[Corollary~8.2]{Mil}.
The general case remains open.
\section{Properties of a Lie group
compared to those of its Lie algebra}
Lie theory derives its strength
from the interplay between
properties
of a Lie group and properties of its Lie algebra.
In the infinite-dimensional case,
the study of
links between $G$ and $L(G)$
has just begun.
It was shown that a connected Lie group~$G$
is abelian if and only if~$L(G)$
is abelian (see \cite[Proposition~22.15]{LEC},
\cite[Proposition~IV.1.10]{MON}
or \cite{GN2}).
Milnor knew this for
regular~$G$.
But for general~$G$
with $L(G)$ abelian,
Milnor stated he could not
prove that $G$ is commutative
\cite[p.\,36]{PRE}.
Neeb achieved essential further progress:
A connected
Lie group~$G$ is solvable (resp., nilpotent)
if and only if $L(G)$ is solvable
(resp., nilpotent)~\cite{GN2}.
\section{Smoothness of continuous homomorphisms}
Milnor asked \cite[p.\,1]{PRE}:
{\em Is a continuous homomorphism
between Lie groups necessarily smooth}\,?
For special types of Lie groups,
this is known:\\[2.5mm]
$\bullet$ Banach-Lie groups (classical);\\[2.5mm]
$\bullet$
Locally exponential
Lie groups~\cite[Lemma~4.3]{PRE};\\[2.5mm]
$\bullet$ Countable direct limits
of finite-dimensional Lie groups
\cite[Prop.\,4.6\,(c)]{FUN}.\\[2.5mm]
Also continuous homomorphisms
from finite-dimensional
Lie groups to
diffeomorphism groups
are smooth
(handwritten notes in \cite{PRE},
also \cite{DIF}; cf.\ \cite[p.\,212]{MaZ}).
Although
Milnor's question remains open,
a positive answer
is available under
stronger
hypotheses:
If a homomorphism $\phi \colon G\to H$ is H\"{o}lder
continuous
at~$1$, then $\phi$ is smooth~\cite[Theorem~3.2]{JFA}.
See \cite[Definition~1.7]{JFA}
for the appropriate concept of H\"{o}lder continuity.
A similar result holds
for the Lie groups of
convenient differential calculus
\cite[Theorem~9.1]{ANN}.
\section{Kernels, Lie subgroups,
quotients and homogeneous spaces}
Milnor asked
whether the kernel of a homomorphism
necessarily is a Lie subgroup
\cite[p.\,1]{PRE},
and proved this for homomorphisms
between locally exponential Lie groups
(they are ``embedded'' Lie subgroups
in the terminology described below).
It is also known
that kernels of smooth homomorphisms
from direct limit groups to Lie groups
are Lie subgroups, like all closed\linebreak
subgroups of such groups~\cite[Proposition~7.5]{FUN}.
But the general answer to Milnor's question
remains open.\\[2.5mm]
We remark that
a wide range of possible concepts of Lie subgroups
is available
in the theory of infinite-dimensional
Lie groups,
each of which can be preferable
in certain situations.
To describe the most basic concept,
let $M$ be a smooth
manifold modelled on a locally convex space~$E$.
A subset $N\subseteq M$
is called a \emph{submanifold}
if there is a
sequentially
closed vector subspace $F\subseteq E$
such that each $x\in N$ is contained in the domain
of some chart
$\phi \colon U\to V\subseteq E$
of~$M$
which takes $U\cap N$ onto $V\cap F$.
Then the restrictions $\phi|_{N\cap U}\colon
U\cap N\to V\cap F$
define a smooth atlas for~$N$.
Note that we do not require that $F$ is complemented
in~$E$ as a topological vector space
(beyond Banach manifolds,
this property loses much of its usefulness).
Given a Lie group~$G$,
a \emph{Lie subgroup} is a subgroup
$H\leq G$ which also is a submanifold.\\[2.5mm]
Also weaker concepts
are needed, analogous
to the ``analytic subgroups''
in finite-dimensional Lie theory.
In the terminology of~\cite{GN2},
an \emph{initial Lie subgroup}
is a subgroup $H\leq G$
which can be given a Lie group structure
which makes the inclusion map
$i\colon H\to G$ a smooth
homomorphism with injective differential
$T_1(i)$, and such that
mappings to~$H$ are smooth if and only if they are smooth
as mappings to~$G$.
It may happen that a subgroup
of a (non-separable)
Banach-Lie group
can be made an analytic subgroup in two
different ways; then one of the Lie group
structures is not the initial one (cf.\ \cite[p.\,157]{HaM}).
Furthermore,
it is not clear whether all subgroups
of interest are initial Lie subgroups.
As a substitute,
one still has the concept of an
\emph{integral subgroup},
referring to an injective
smooth homomorphism
$i\colon H\to G$ from a (connected)
Lie group to~$G$ such that $T_1(i)$
is injective.
Milnor used the notion of an
\emph{immersed Lie subgroup}:
this is an injective smooth
homomorphism of Lie groups
$i\colon H\to G$
taking some open identity
neighbourhood in~$H$ onto a
submanifold of~$G$ (cf.\ \cite[p.\,22]{PRE}).\\[2.5mm]
Also stronger notions of Lie subgroups
are needed.
A Lie subgroup $H\leq G$ is called
a \emph{split} if
$G/H$ can be given a smooth
manifold structure making
the canonical map $q\colon G\to G/H$
a smooth $H$-principal bundle
(i.e., $q$ is smooth and
admits smooth local sections).
For $G$ locally exponential,
the concept of an \emph{embedded
Lie subgroup}~$H$
is particularly useful.
This is a sequentially closed
subgroup such that $H\cap\exp_G(U)
=\exp_G(U\cap {\mathfrak h})$
for a $0$-neighbourhood $U\subseteq L(G)$
on which $\exp_G$ is injective,
where
\[
{\mathfrak h}\; :=\; \{X\in L(G)\colon
\exp_G({\mathbb R} X)\subseteq H\}\, .
\]
For example, it can be shown that the topological
quotient group $G/N$ of a BCH-Lie group~$G$
modulo a closed normal subgroup~$N$ of~$G$
is a BCH-Lie group if and only if
$N$ is an embedded Lie subgroup of~$G$
(\cite[Corollary~2.21]{GCX};
cf.\ \cite{GN1} for Banach-Lie groups).
This result was extended to
locally exponential Lie groups~$G$
by Neeb~\cite{GN2}.
In this case, $G/N$ is a locally
exponential Lie group if and only if
$N$ is an embedded Lie subgroup
whose Lie algebra $L(N)$ is ``locally exponential.''
Neeb also showed that
every locally compact subgroup
of a locally exponential Lie group
is an embedded Lie subgroup~\cite{GN2},
as in the case of Banach-Lie groups
(first discussed by Birkhoff).\\[2.5mm]
It would be very useful to find tangible criteria
ensuring that a homogenous space $G/H$
can be made a smooth manifold with
reasonable properties,\footnote{A minimal requirement is that
the smooth manifold structure on $G/H$ is final
with respect to the quotient
map $q\colon G\to G/H$.
In addition to this, one
would certainly like to require
that $T_1(q)$ is
a quotient homomorphism with kernel $T_1(H)$.}
for a closed subgroup
$H$ of a locally
exponential Lie group (which does not happen to
be normal or a split Lie subgroup).
Such criteria are not even known in the case
of Banach-Lie groups, and any progress
in this direction would be most valuable.\\[2.5mm]
Also, it would be desirable
to clarify the precise relations between the various
concepts of Lie subgroups,
and to find examples which clearly distinguish
the concepts.
For instance,
it is well known that embedded Lie subgroups
and ordinary Lie subgroups
of Banach-Lie groups coincide.
But it is unclear whether
these concepts still agree in the case of
BCH- (or locally exponential) Lie groups.
\section{Integrability questions}
Milnor \cite[p.\,1]{PRE}
asked whether every closed subalgebra
of $L(G)$ corresponds to some
immersed Lie subgroup of~$G$.
In \cite[Warning~8.5]{Mil}, he
described a counterexample
(due to Omori).
The integrability
question of Lie subalgebras
was analyzed further in~\cite{Rob}
and \cite{GN2},
notably for locally exponential Lie
groups.\\[2.5mm]
It is a classical result by
van Est and Korthagen
that a Banach-Lie algebra ${\mathfrak g}$
need not be ``integrable'' (or ``enlargible'')
-- there need
not be a (Banach-) Lie group~$G$
with $L(G)\cong {\mathfrak g}$.
As they showed, ${\mathfrak g}$ is integrable
if and only if a certain subgroup
$\Pi({\mathfrak g})\leq {\mathfrak z}({\mathfrak g})$
of the center of ${\mathfrak g}$
(the period group) is discrete~\cite{EaK}.
Related to this work is a
question by Milnor \cite[pp.\,31-32]{PRE},
which can be re-phrased
as follows:
If $G$ is a BCH-Lie group,
does it follow that $L(G)_{\mathbb C}$
is integrable to a Lie group\,?
The answer is negative,
even for Banach-Lie groups
(see \cite[Example~VI.4]{GN1}).\\[2.5mm]
Milnor remarks that it would be interesting
to know which topological Lie algebras~${\mathfrak g}$
correspond to BCH-Lie groups.
A necessary condition is
that the BCH-series
converges on a $0$-neighbourhood
of ${\mathfrak g}\times {\mathfrak g}$
to an analytic function
(see~\cite{Rb2} for characterizations
of this property).
A full solution to Milnor's
question was given by Neeb,
even for the wider class
of ``locally exponential'' Lie algebras.
Any such Lie algebra~${\mathfrak g}$ can be associated
a certain period group $\Pi({\mathfrak g})
\leq {\mathfrak z}({\mathfrak g})$;
it is integrable to a locally exponential Lie group
if and only if $\Pi({\mathfrak g})$
is discrete~\cite{GN3}.\\[2.5mm]
Integrability questions
have also been studied
for other types of Lie algebras.
It was shown that
every locally finite Lie algebra
of countable dimension
is integrable~\cite[Theorem~5.1]{FUN}.
In \cite{Nee}, Neeb described
the obstructions
to integrate a central extension
of topological Lie algebras
to a Lie group extension.
Later, he
extended his methods to abelian~\cite{Ne2}
and non-abelian extensions~\cite{Ne3}.
{\footnotesize
{\bf Helge Gl\"{o}ckner}\\[.4mm]
Darmstadt University of Technology\\
Department of Mathematics, AG~5,\\
Schlossgartenstr.\ 7\\
64289 Darmstadt, Germany\\[1.5mm]
{\it E-mail address}:\\
\,{\tt gloeckner@mathematik.tu-darmstadt.de}}
\end{document}